\documentclass{llncs}
\usepackage{llncsdoc}
\usepackage{epsfig}
\usepackage{graphicx}
\usepackage[leqno]{amsmath}
\usepackage{amssymb}
\usepackage{mathrsfs}
\usepackage{psfrag}
\usepackage[usenames]{color}
\newcommand{\R}{\mathbb R}
\newcommand{\Z}{\mathbb Z}
\newcommand{\C}{\mathbb C}
\newcommand{\eps}{\varepsilon}

\begin{document}
\title{Comparing persistence diagrams through complex vectors}

%
\author{Barbara Di Fabio\inst{1} \and Massimo Ferri\inst{1,2}}
%
%
%
\institute{ARCES, University of Bologna, Italy\\
\email{barbara.difabio@unibo.it},
\and
Department of Mathematics, University of Bologna, Italy\\
\email{massimo.ferri@unibo.it}}

\maketitle              

\begin{abstract}
The natural pseudo-distance of spaces endowed with filtering functions is precious for shape classification and retrieval; its optimal estimate coming from persistence diagrams is the bottleneck distance, which unfortunately suffers from combinatorial explosion. A possible algebraic representation of persistence diagrams is offered by complex polynomials; since far polynomials represent far persistence diagrams, a fast comparison of the coefficient vectors can reduce the size of the database to be classified by the bottleneck distance. This article explores experimentally three transformations from diagrams to polynomials and three distances between the complex vectors of coefficients.

\end{abstract}
\section{Introduction}\label{intro}

Persistent homology has already proven to be an effective tool for
shape representation in various applications, in particular when
the objects to be classified, compared or retrieved have a natural
origin. The interplay of geometry and topology in persistence
makes it possible to capture qualitative aspects in a formal and
computable way, yet it doesn't suffer of the excessive freedom of
mere topological equivalence. The privileged tool for shape
comparison is the natural pseudo-distance \cite{DoFr04}, which is scarcely
computable. Luckily, persistence diagrams condense the essence of
the shape concept of the observer in finite sets of points in the
plane \cite{FrLa01,EdHa08}; moreover, the bottleneck distance (a.k.a. matching distance) between persistence
diagrams yields an optimal lower bound to the natural
pseudo-distance \cite{dAFrLa,DoFr04bis}. There is a problem: the  bottleneck distance suffers
from combinatorial explosion \cite{dAFrLabis}, so it becomes hard to scan a large
database when it comes to retrieval. Approximations, smart
organization of the database according to the metric, progressive
application of different classifiers come to help, but the problem
is lightened, not solved.

A paradigm shift came from an idea of Claudia Landi \cite{FeLa99}:
represent the persistence diagram as the set of complex roots of a
polynomial; then comparison can be performed on coefficients. Two
problems arise: one --- which comes from the nature itself of persistence diagrams ---
is that in real situations there are a lot of
points near the ``diagonal'' $\{(u,v)\in\R^2: u=v\}$, due to noise so less
meaningful in shape representation, but with a heavy impact on polynomial coefficients;
another problem --- coming from polynomial theory --- is that little
distance of polynomial roots implies little distance of
coefficients, but the converse is false.

A completely different polynomial representation of barcodes
(equivalently: of persistence diagrams) is the one through
tropical algebra \cite{KaCa15}, closely adapting to the bottleneck
distance. \bigskip

\noindent\textbf{The contribution of the paper.} We face the first
problem --- the existence of points near the diagonal --- by
performing a plane warping which takes all the line $u=v$ to 0, so
points near the diagonal actually become close together. Making
noise points close and around zero diminishes their contribution
to polynomial coefficients, above all to the first (and most
relevant) ones: sum of roots, sum of pairwise products of roots,
etc. As for the second problem --- the fact that close
coefficients may not mean close roots --- we explore the use of
polynomial comparison as a preprocessing phase in shape retrieval,
i.e. as a very fast way of getting rid of definitely far objects,
so that the bottleneck distance can be computed only for a small
set of candidates, in the same line as \cite{CeDiJaMe14}. The
results are satisfactory: in some of our experiments the
bottleneck distance even turns out to be unnecessary.

\section{Preliminaries}

In persistence, the shape of an object is usually studied by
choosing a topological space $X$ to represent it, and a function
$f:X\to \R$, called a {\em filtering} (or {\em measuring}) {\em
function}, to define a family of  subspaces
$X_u=f^{-1}((-\infty,u])$, $u\in\R$, nested by inclusion, i.e. a
filtration of $X$. The name ``persistence'' is bound to the idea
of ranking topological features by importance, according to the
length of their ``life'' through the filtration. The basic
assumption is that the longer a feature survives, the more
meaningful or coarse the feature is for shape description. In
particular, structural properties of the space $X$ are identified
by features that once born never die; vice-versa, noise and shape
details are characterized by a short life. To study how
topological features vary in passing from a set of the filtration
into a larger one we use homology. A nice feature of this approach
is modularity: The choice of different filtering functions may
account for different viewpoints on the same problem (different
shape concepts) or for different tasks. For further details we
refer to \cite{BiDe*08,EdHa09}.

Persistent homology groups of the pair $(X, f)$ --- i.e. of the
filtration $\{X_u\}_{u\in\R}$ --- are defined as follows. Given
$u\leq v\in\R$, we consider the inclusion of $X_u$ into $X_v$.
This inclusion induces a homomorphism of homology groups
$H_k(X_u)\rightarrow H_k(X_v)$ for every $k\in\Z$. Its image
consists of the $k$-homology classes that are born before or at
the level $u$  and are still alive at the level $v$ and is called
the {\em $k$th persistent homology group of $(X,f)$ at $(u, v)$}.
When this group is finitely generated, we denote by
$\beta_k^{u,v}(X,f)$ its rank.

The usual, compact description of persistent homology groups of
$(X,f)$ is provided by the so-called \emph{persistence diagrams},
i.e. multisets of points whose abscissa and ordinate are,
respectively, the level at which $k$-homology classes are created
and the level at which they are annihilated through the
filtration. If a homology class does not die along the filtration,
the ordinate of the corresponding point is set to $+\infty$.

At the moment, our approach to convert persistence diagrams into
complex vectors can be applied only when neglecting these points
at infinity. Hence, we focus on the subsets of proper points of
the classical persistence diagrams, known in literature as
\emph{ordinary persistence diagrams} \cite{CoEdHa09}. For
simplicity we still call them ``persistence diagrams''. We
underline that this choice is not so restrictive since the number
of points at infinity depend only on the homology of the space
$X$, and persistent homology provide a finite distance between two
pairs if and only if the considered spaces are homeomorphic.

We use the following notation: $\Delta^+=\{(u,v)\in \R^2: u<v\}$,
$\Delta=\{(u,v)\in \R^2: u=v\}$, and
$\overline{\Delta^+}=\Delta^+\cup\Delta$.

\begin{definition}
\label{cornerpt} Let $k\in\Z$ and $(u,v)\in\Delta^+$. The
\emph{multiplicity} $\mu_k (u,v)$  of $(u,v)$ is the
finite non-negative number defined by {\setlength\arraycolsep{1pt}
\begin{equation*}
\lim_{\eps\to 0^+}\left(\beta_k^{u+\eps
,v-\eps}(X,f)-\beta_k^{u-\eps ,v-\eps}(X,f)-\beta_k^{u+\eps
,v+\eps}(X,f)+\beta_k^{u-\eps ,v+\eps}(X,f)\right).
\end{equation*}}
\end{definition}
\begin{definition}\label{persdiag}
The \emph{$k$th-persistence diagram} $D_k(X, f)$ is the set of all points
$(u,v)\in\Delta^+$ such that $\mu_k(u,v)>0$, counted with
their multiplicity, union the points of $\Delta$, counted with
infinite multiplicity. We call \emph{proper points} the points of a persistence diagram lying on $\Delta^+$.
\end{definition}

\begin{figure}[ht]
\psfrag{D}{$D_0(X,f)$}\psfrag{a}{$a$}\psfrag{b}{$b$}\psfrag{c}{$c$}\psfrag{d}{$d$}
\psfrag{x}{$u$}\psfrag{y}{$v$}\psfrag{p}{$p_1$}\psfrag{q}{$p_3$}\psfrag{r}{$p_2$}
\psfrag{X}{$X$}\psfrag{f}{$f$}
\begin{center}
\includegraphics[width=8cm]{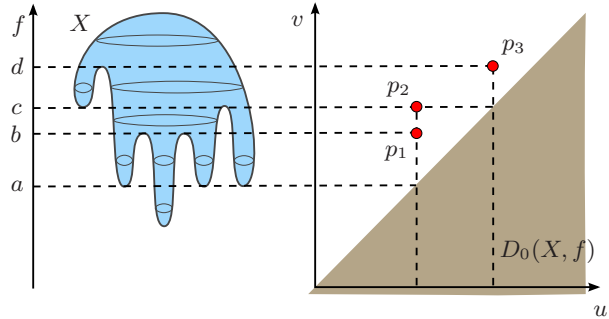}
\end{center}
\caption{\footnotesize{Left: The height function $f$ on the space
$X$. Right: The associated $0th$-persistence diagram
$D_0(X,f)$.}}\label{PersDiagExample}
\end{figure}

Figure~\ref{PersDiagExample} displays an example of persistence
diagram for $k=0$. The surface $X\subset\R^3$ is filtered by
the height function $f$. $D_0(X,f)$ has three proper points $p_1, p_2, p_3$ since the abscissa of these points
corresponds to the level at which new connected components
are born along the filtration, while the
ordinate identifies the level at which these connected components
merge with existing ones. In terms of multiplicity, this means that $\mu_0(p_i)>0$, $i=1,2,3$, and $\mu_0(p)=0$ for every other point $p\in\Delta^+$. To see, for example, that
$\mu_{0}(p_1)=2$, where $p_1=(a,b)$, it is sufficient to
observe that, for every $\eps>0$ sufficiently small, it holds that
$\beta_0^{a+\eps ,b-\eps}(X,f)=4$, $\beta_0^{a+\eps ,b+\eps}(X,f)=2$, $\beta_0^{a-\eps ,b-\eps}(X,f)=\beta_0^{a-\eps ,b+\eps}(X,f)=1$, and
apply Definition \ref{cornerpt}. In an analogous way, it can be
observed that $\mu_0(p_2)=\mu_0(p_3)=1$.

Persistence diagrams comparison is usually carried out through the so called \emph{bottleneck distance} because of the robustness of these descriptors with
respect to it. Roughly, small changing in a given
filtering function (w.r.t. the max-norm) produces just a small
changing in the associated persistence diagram w.r.t. the bottleneck
distance \cite{CoEdHa07,CeDiFeFrLa09}.
The bottleneck distance between two persistence diagrams measures
the cost of finding a correspondence between their points. In
doing this, the cost of taking a point $p$ to a point $p'$ is
measured as the minimum between the cost of moving one point onto
the other and the cost of moving both points onto the diagonal. In particular,
the matching of a proper point $p$ with a point of $\Delta$ can be
interpreted as the destruction of the point $p$. Formally:
\begin{definition}\label{matchDist}
Let $D$, $D'$ be two persistence diagrams. The bottleneck
distance $d_B\left(D,D'\right)$ is defined as
\begin{equation*}
d_B(D,D')=\min_{\sigma}\max_{p\in
D}d(p,\sigma(p)),
\end{equation*}
where $\sigma$ varies among all the bijections between $D$ and
$D'$ and
$$\label{deltaDistance}
d\left(\left(u,v\right),\left(u',v'\right)\right)=\min\left\{\max\left\{|u-u'|,|v-v'|\right\},\max\left\{\frac{v-u}{2},\frac{v'-u'}{2}\right\}\right\}
$$
for every
$\left(u,v\right),\left(u',v'\right)\in\overline{\Delta^+}$.
\end{definition}

\section{Persistence diagrams vs complex vectors}\label{vectors}
Driven by the awareness that, in the experimental framework,
evaluating the bottleneck distance can be computationally
expensive, making its usage not practicable on large datasets, in
this work we propose a new procedure based on the preliminary idea
introduced in \cite{FeLa99}. We translate the problem of comparing
directly two persistence diagrams through the bottleneck distance
into the problem of comparing complex vectors associated with each
persistence diagram through appropriate metrics between vectors.
The components of these complex vectors are complex polynomials'
coefficients obtained as follows. Firstly, we define a certain
transformation taking points of persistence diagrams to the set of
complex numbers. Secondly, we construct a complex polynomial
having the obtained complex numbers as roots.

In this paper, we consider the three transformations below:

\begin{itemize}
\item  $R:\overline{\Delta^+}\to\C$, with $R(u,v)=u+iv$,
\item  $S:\overline{\Delta^+}\to\C$, with $S(u,v)=\left\{\begin{array}{ll}\dfrac{v-u}{\alpha\sqrt{2}}\cdot(u+iv),& \mbox{if}\, (u,v)\neq (0,0)\\(0,0),&\mbox{otherwise}\end{array}\right.$,
\item $T:\overline{\Delta^+}\to\C$, with $T(u,v)=\dfrac{v-u}{2}\cdot(\cos\alpha -\sin\alpha+i(\cos\alpha+\sin\alpha))$,
\end{itemize}
where $\alpha=\sqrt{u^2+v^2}$.

\begin{figure}[ht]
\psfrag{S(P)}{\!\!$S(p)$}\psfrag{P}{\!$p$}\psfrag{T(P)}{\!\!$T(p)$}\psfrag{0}{$0$}\psfrag{2}{$2$}\psfrag{-2}{\!\!\!$-2$}
\psfrag{4}{$4$}\psfrag{y}{$v$}\psfrag{e}{$\varepsilon$}\psfrag{d}{$\delta$}\psfrag{p}{$p$}
\begin{center}
\includegraphics[width=\textwidth]{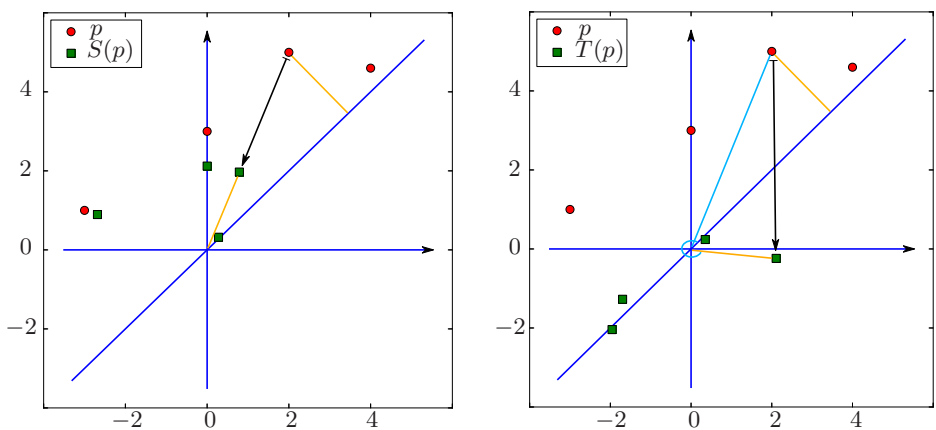}
\end{center}
\caption{A persistence diagram with its transformations $S$ (left)
and $T$ (right). Same colors identify same
lengths.}\label{ExampleST}
\end{figure}

$R,S,T$ are continuous maps; $R$ and $S$ are also injective on
$\overline{\Delta^+}$ and $\Delta^+$, respectively. We define the
multiplicity of a complex number in the range of $R,S,T$ to be the
sum of the multiplicities of the points belonging to its preimage
(this is necessary because of the non-injectivity of $T$ on
$\Delta^+$, although a preimage containing more than one proper
point of the diagram has zero probability to occur). The main
differences among these deformations are the following: the
deformation $R$ acts as the identity, just passing from $\R^2$ to
$\C$; the deformation $S$ warps the diagonal $\Delta$ to the
origin, and takes points of $\Delta^+$ to points of  $\{z\in\C:
Re(z)<Im(z)\}$; the deformation $T$ warps the diagonal $\Delta$ to
the origin, and takes points of $\Delta^+$ to points of $\C$. An
example showing how $S$ and $T$ transform a persistence diagram is
represented in Figure \ref{ExampleST}.
Both $S$ and $T$ seemed to be preferable to $R$ because points
near $\Delta$ --- due to noise --- have to be considered close to
each other in the bottleneck distance, although they may be very
far apart in Euclidean distance. Taking them to be all near the
origin would then also reduce their impact in the sums and sums of
products which will build the polynomial coefficients we are going
to compare. In particular, $T$ was designed to distribute the
image of those noise points around zero, whereas $S$ makes them
near zero, but all on one side: in the half-plane of $\C$
corresponding to $\Delta^+$. $T$ has two drawbacks: it is not
injective on $\Delta^+$ and does not behave well with respect to
simple transformations.

Let $D$ be a persistence diagram, and $p_1=(u_1, v_1), \ldots,
p_s=(u_s, v_s)$ its proper points with multiplicity $r_1, \ldots,
r_s$, respectively. Let now the complex numbers $z_1, \ldots, z_s$
be obtained from $p_1, \ldots, p_s$ by one of the transformations
$R$, $S$ or $T$. We associate to $D$ the complex polynomial
$f_D(t) = \prod_{j=1}^s (t-z_j)^{r_j}$. We are actually interested
in the coefficient sequence of $f_D(t)$, which we can compute by
Vi\`ete's formulas (see Algorithm 2).

Once we have the polynomials $f_{D}(t) = t^n-a_1 t^{n-1} + \cdots
+ (-1)^i a_i t^{n-i} + \cdots + (-1)^n a_n$ and $f_{D'}(t) =
t^m-a'_1 t^{m-1} + \cdots + (-1)^j a'_j t^{m-j} + \cdots + (-1)^m
a'_m$ corresponding to persistence diagrams $D$, $D'$, we face a
first problem, given by the possibly different degrees $n$ and $m$
($m<n$ say). Because of their expression in terms of roots, we
prefer to compare coefficients with the same index, rather than
coefficients relative to the same degree of $t$. We manage this
problem by adding $n-m$ null coefficients to $f_{D'}(t)$, i.e.
multiplying $f_{D'}(t)$ by $t^{n-m}$, which amounts to adding the
complex number zero with multiplicity $n-m$. In so doing, we can
build two vectors of complex numbers
$(a_1,\dots,a_n),(a'_1,\ldots,a'_n)$ of the same length and are
ready to compute a distance between them. By continuity of
Vi\`ete's formulas, close roots imply close coefficients. Hence,
two persistence diagrams that are close in terms of bottleneck
distance have close associated polynomials. Unfortunately, the
converse is not true.

Preliminary tests suggested that the first coefficients were more
meaningful; therefore we experimented with different distances on
two complex vectors $(a_1,\dots,a_k),(a'_1,\ldots,a'_k)$,
$k\in\{1,\ldots,n\}$, one treating all coefficients equally, two
which give decreasing value to coefficients of increasing indices.
The chosen metrics are the following:
$$d_1=\underset{j=1}{\overset{k}\sum}|a_j-a'_j|,\quad d_2=\underset{j=1}{\overset{k}\sum}\dfrac{|a_j-a'_j|}{j},\quad d_3=\underset{j=1}{\overset{k}\sum}|a_j-a'_j|^{1/j}.$$

\noindent\textbf{Algorithms and computational analysis.} The algorithms below resume the principal steps of our scheme. $F$ in Algorithm 1 (line 2) and  $d$ in Algorithm 3 (line 4) correspond, respectively, to one of the transformations $R,S,T$ and one of the metrics $d_1,d_2,d_3$ previously defined.

\begin{tabular}{l}
\textbf{Algorithm 1:} ComplexLists\\
\textbf{Input:} List $A$ of proper points of a persistence diagram $D$,\\
$M=\max_{|A|}\{A: A\in \,\mbox{database}\, Db\}$\\
\textbf{Output:} List $B$ of complex numbers associated with $D$\\
\begin{tabular}{r l|r l}
\hline
1:& \  \textbf{for each} $(u,v)\in A$&4:& \  \textbf{if} $|B|<M$\\
2:& \   \quad \textbf{replace} $(u,v)$ by $F(u,v)$&5:& \   \quad \textbf{append} $M-|B|$ zeros to $B$\\
3:& \  \textbf{end for}&6:& \  \textbf{end if}\\
\hline\\
\end{tabular}\\
\end{tabular}\label{CL-Alg}

\begin{tabular}{l}
\textbf{Algorithm 2:} ComplexVectors\\
\textbf{Input:} $M$, $B=list(z_1,\ldots,z_M)$ associated with $D$, $k\in [0,M]$\\
\textbf{Output:} Complex vector $V_k$ associated with $D$\\
\begin{tabular}{r l}
\hline
1:& \  \textbf{set} $V_k=list()$\\
2:& \   \textbf{for} $j\in \{1,\ldots,k\}$\\
3:& \   \quad \textbf{compute} $c_j(z_1,\ldots,z_M)=\underset{1\le i_1<i_2<\ldots<i_j\le M}\sum z_{i_1}\cdot z_{i_2}\cdot\ldots\cdot z_{i_j}$\\
4:& \  \quad \textbf{append} $c_j$ to $V_k$\\
5:& \   \textbf{end for}\\
\hline\\
\end{tabular}
\end{tabular}\label{CV-Alg}

\begin{tabular}{l}
\textbf{Algorithm 3:} VectorsComparison\\
\textbf{Input:} $L=\{V_k: V_k \,\mbox{complex vector associated with}\, D \,\mbox{for each}\, D\in Db\}$\\
\textbf{Output:} Matrix of distances $d(V_k,V'_k)$\\
\begin{tabular}{r l|r l}
\hline
1:& \  \textbf{set} $M=(0_{ij})$, $i,j={1,\ldots,|L|}$&4:& \   \quad\quad \textbf{replace} $0_{ij}, 0_{ji}$ by $d(i,j)$\\
2:& \  \textbf{for each} $i\in \{1,\ldots,|L|\}$&5:& \  \quad \textbf{end for}\\
3:& \   \quad\textbf{for each} $j\in\{i,\ldots,|L|\}$&6:& \   \textbf{end for}\\
\hline\\
\end{tabular}
\end{tabular}\label{VC-Alg}

Let $N=|L|=|Db|$. It is easily seen that the computational
complexities of Algorithms 1 and 3 are $C_1=O(M\cdot N)$ and
$C_3=O(k\cdot N^2)$, respectively. The cost of Algorithm 2 depends
on how we have implemented the computation of Vi\`ete formulas.
Using the induction on the index $j$, we have $C_2= O((2k^2+k\cdot M)\cdot N)$.

We want to show that our approach to database classification, in
general, results to be cheaper than using the bottleneck distance
by a suitable choice of the number $k$ of computed coefficients.
Our comparison is realized in terms of storage locations and not
in terms of running time performances since the algorithm proposed
here and the one based on the bottleneck distance run on different
platforms. We recall that the cost of computing the bottleneck
distance $d_B(D,D')$ is $C_B=O\left((r+r')^{3/2}\log(r+r')\right)$
if $A,A'$ are the subsets of proper points of two persistence
diagrams $D,D'$ with $|A|=r, |A'|=r'$ (see \cite{EfItKa01}).
Instead, using our scheme, with $N=2$ and $M=\max(r,r')$, we get
$O((\max(r,r')+2k^2+k\cdot\max(r,r')+2k)$, so
$C=O(k\cdot(\max(r,r')+k))$. Since $k\le\max(r,r')$, in the worst
case, we have $C=O\left((\max(r,r'))^2\right)$ which is higher
than $C_B$, but for pre-processing we may choose a favorable $k$
(e.g. $k=\left\lfloor \sqrt{\max(r,r')}\right\rfloor$). Also
consider that, for a retrieval task, the heavy part of the
computation (Vi\`ete's formulas) for the database is performed
offline; in other words: if we store the coefficient vectors
instead of the proper points lists, then the search can be
performed by a distance computation of complexity
$O\left(\max(r,r')\right)$!

\section{Experimental results}\label{Exp}
This section is devoted to validate the theoretical framework
introduced in Section \ref{vectors}. In particular, through some
experiments on persistence diagrams for 0th homology degree
associated with 3D-models represented by triangle meshes, we will
prove that our approach allows to perform the persistence diagrams
comparison without greatly affecting (and in some cases improving)
the goodness of results in terms of database classification.

\begin{figure}[htbp]
\psfrag{w}{$w$}\psfrag{B}{\textcolor[rgb]{1.00,0.00,0.00}{$B$}}\psfrag{P}{\textcolor[rgb]{0.00,0.67,0.00}{$P$}}\psfrag{L}{\textcolor[rgb]{0.00,0.00,1.00}{$L$}}
\begin{center}
\includegraphics[width=9cm]{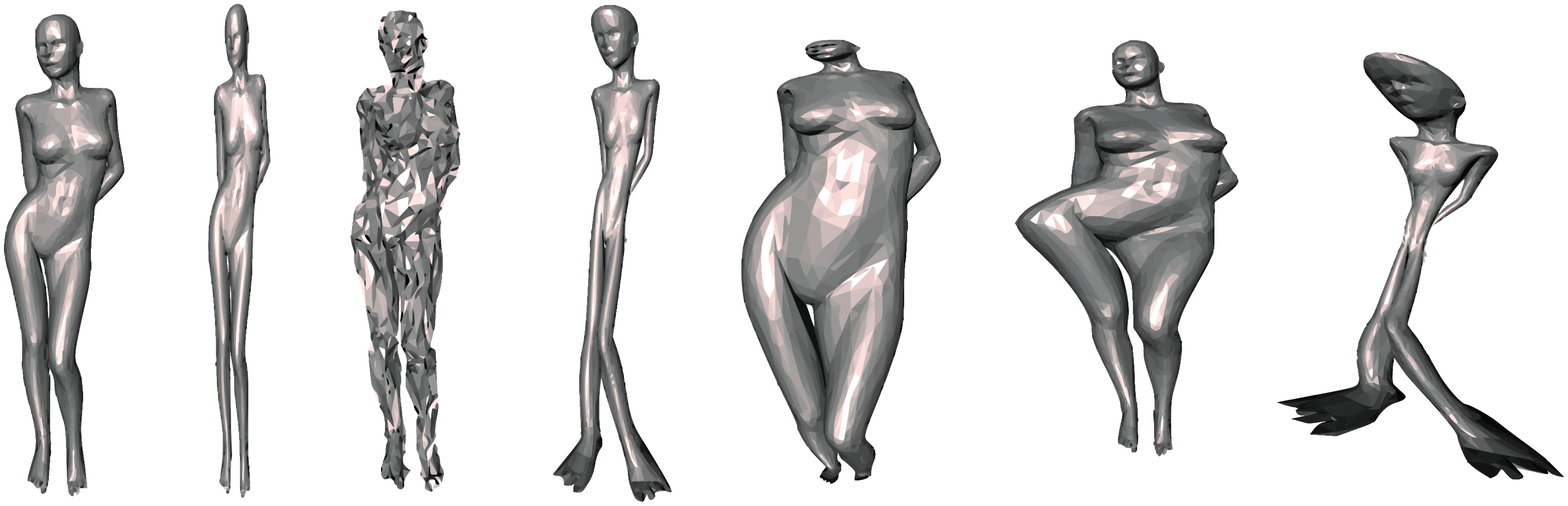}\\
\caption{\footnotesize{The null model ``Victoria0'' and the 3rd strength level for each deformation.}}\label{victoria}
\end{center}
\begin{center}
\includegraphics[width=5cm]{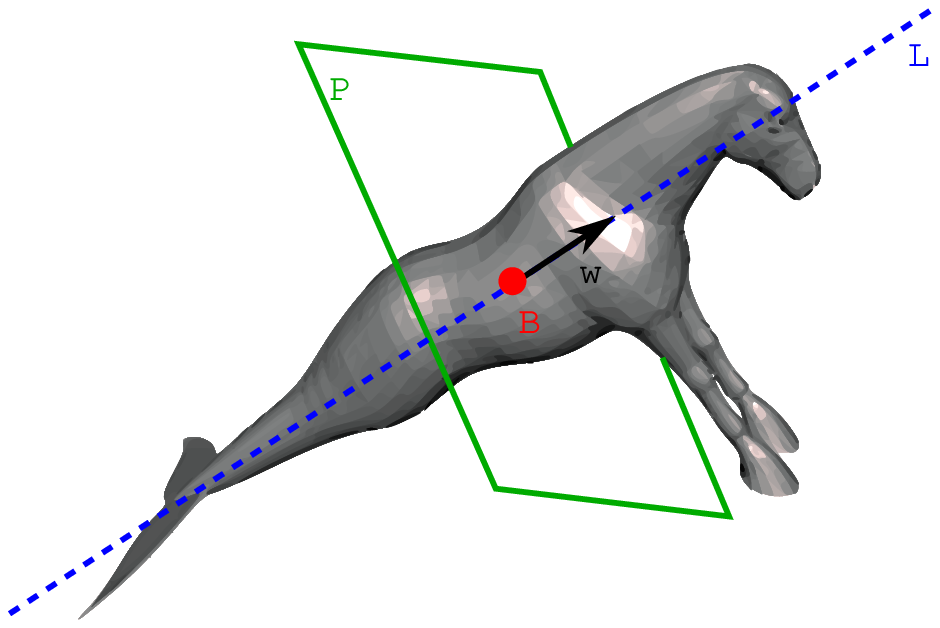}\\
\caption{\footnotesize{The null model ``seahorse0'' depicted with its center of mass $B$ and the associated vector $w$,
which define the filtering functions $f_L$ and $f_P$.}}\label{seahorse}
\end{center}
\end{figure}

To test the proposed framework we considered a database of 228
3D-surface mesh models introduced in \cite{BiCe*11}. The database
is divided into 12 classes, each containing 19 elements obtained
as follows: A null model taken from the Non Rigid World Benchmark
\cite{linkTosca} is considered together with six non-rigid
transformations applied to it at three different strength levels.
An example of the transformations and their greatest strength levels is
given in Figure~\ref{victoria}.

Two filtering functions $f_L$, $f_P$ have been
defined on the models of the database as follows: For each triangle mesh $M$ of vertices $\{v_1,\dots,v_n\}$,
the center of mass $B$ is computed, and the model is normalized to
be contained in a unit sphere. Further, a vector $w$ is
defined as
\begin{equation*}
w=\frac{\sum_{i=1}^n(v_i-B)\|v_i-B\|}{\sum_{i=1}^n\|v_i-B\|^2}.
\end{equation*}
The function $f_L$ is the distance from
the line $L$ parallel to $w$ and passing through $B$, while the function $f_P$
is the distance from the plane $P$ orthogonal to $w$ and passing
through $B$ (see, as an example, Figure~\ref{seahorse}). The values
of $f_L$ and $f_P$ are then normalized so
that they range in the interval $[0,1]$. These filtering functions
are translation and rotation invariant, as well as scale invariant
because of a priori normalization of the models. Moreover, the
considered models are sufficiently generic (no point-symmetries
occur etc...) to ensure that the vector $w$ is well-defined
over the whole database, as well as its orientation stability.

Our experimental results are synthesized in Tables \ref{graphL}
and \ref{graphP} in terms of precision/recall (PR) graphs when the
filtering functions $f_L, f_P$, respectively, are considered.
Before going into details, we want to emphasize that our intent is
not to validate the usage of persistence for shape comparison,
retrieval or classification. In fact, as a reader coming from the
retrieval domain will probably note, the PR graphs reported in
this paper are below the state of the art. This depends on the
fact that, in general, good retrieval performances can be achieved
only taking into account different filtering functions that give
rise to a battery of descriptors associated with each model in the
database.

\begin{table}[htbp]
\begin{center}
\caption{\footnotesize{PR graphs related to the filtering function $f_L$, when $0th$-persistence diagrams are compared directly through the bottleneck distance and in terms of the first $k$ components of the complex vectors obtained from the transformations $R$ (first row), $S$ (second row) and $T$ (third row) through the distances $d_1$ (first column), $d_2$ (second column) and $d_3$ (third column).}}\label{graphL}
\begin{tabular}{ccc}
\includegraphics[height=3.2cm]{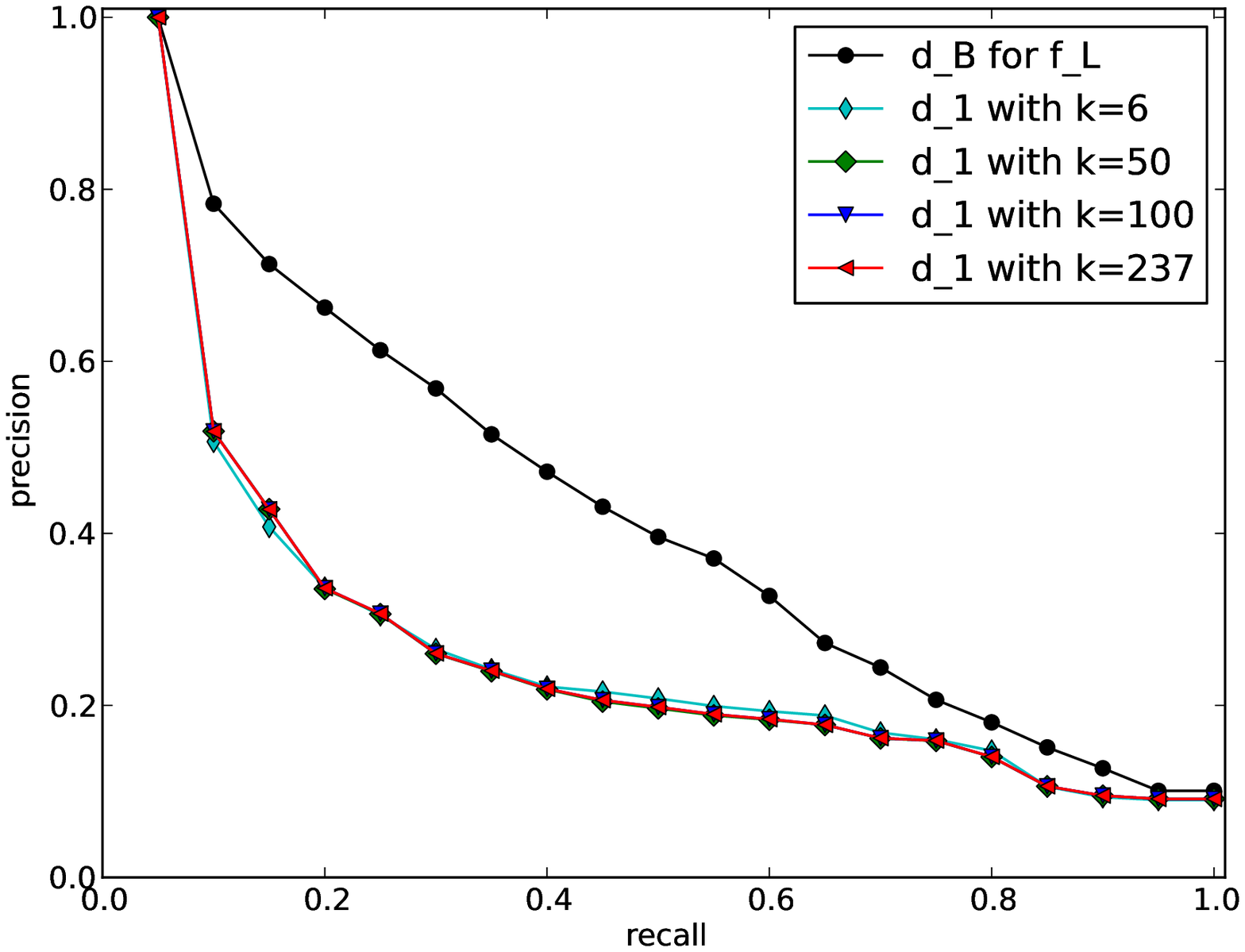}\!\!\!\!&\!\!\!
\includegraphics[height=3.2cm]{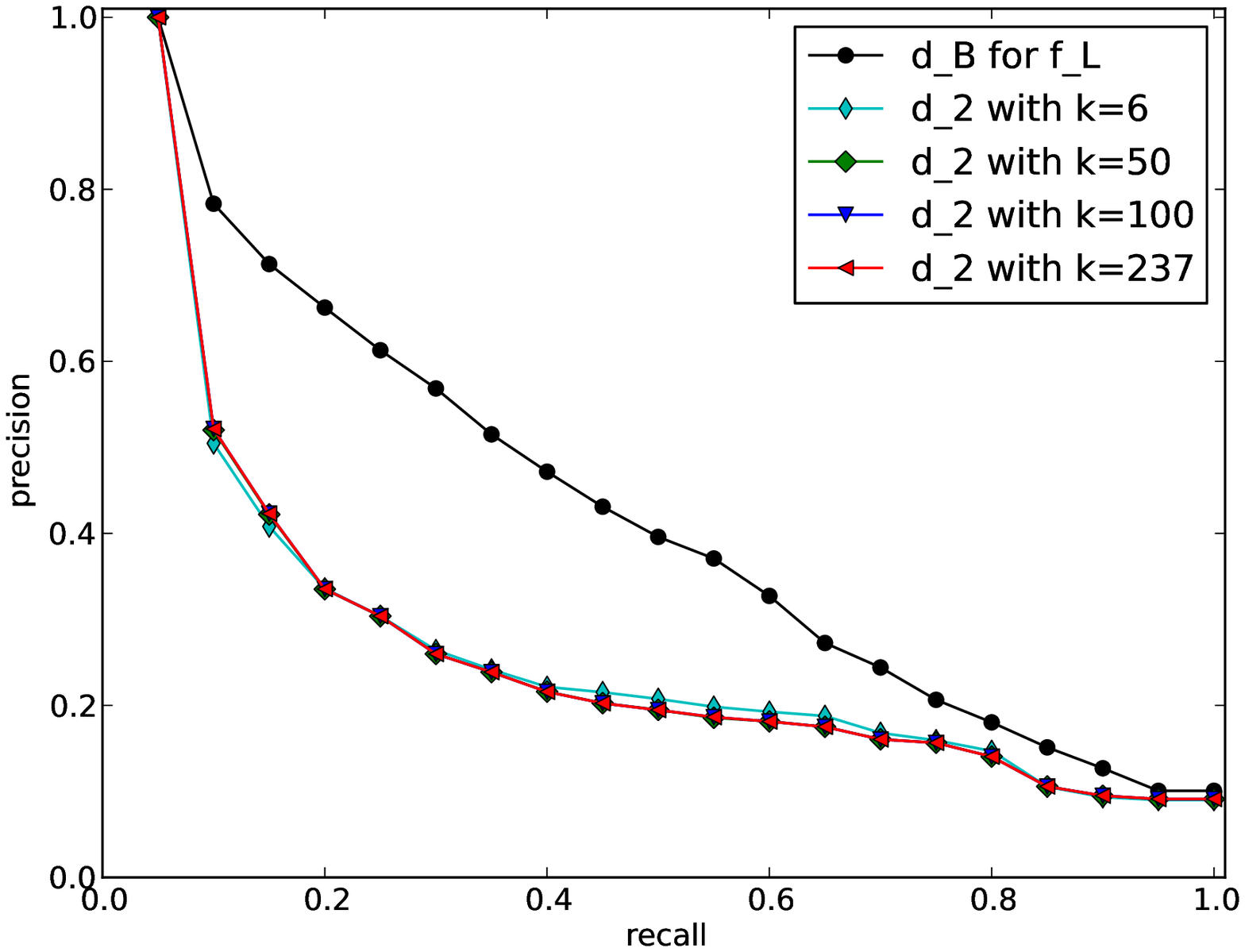}\!\!\!\!&\!\!\!
\includegraphics[height=3.2cm]{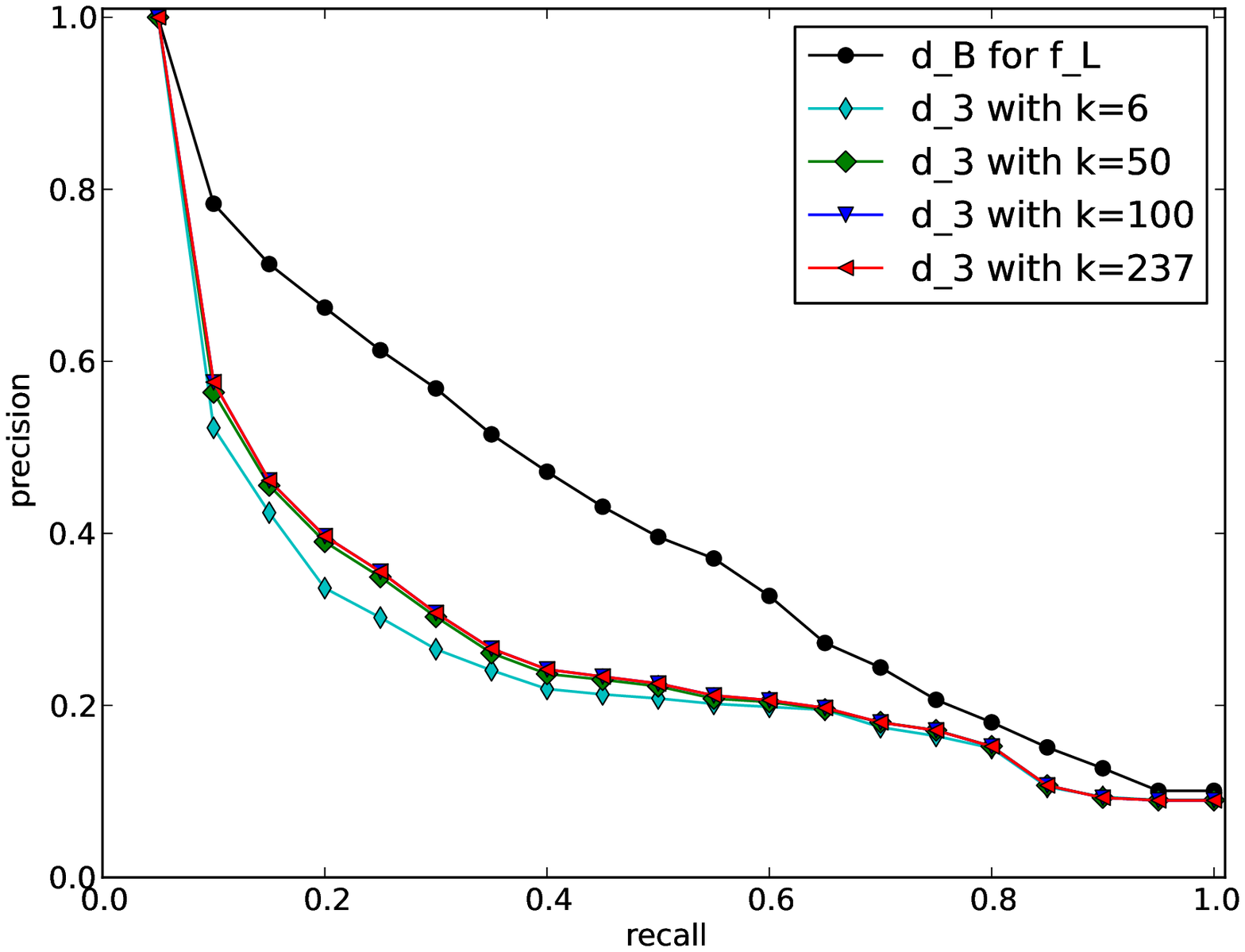}\\
\includegraphics[height=3.2cm]{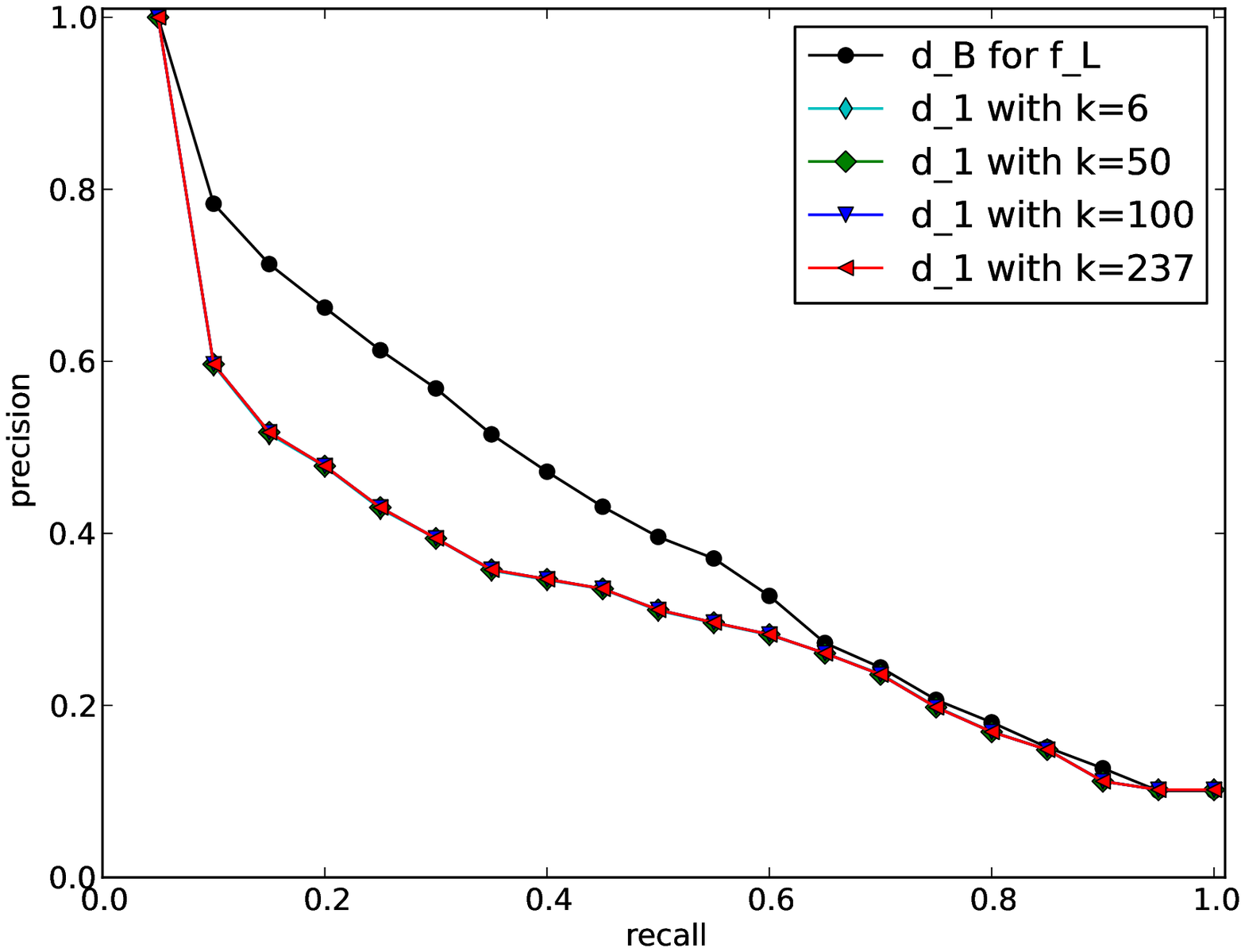}\!\!\!\!&\!\!\!
\includegraphics[height=3.2cm]{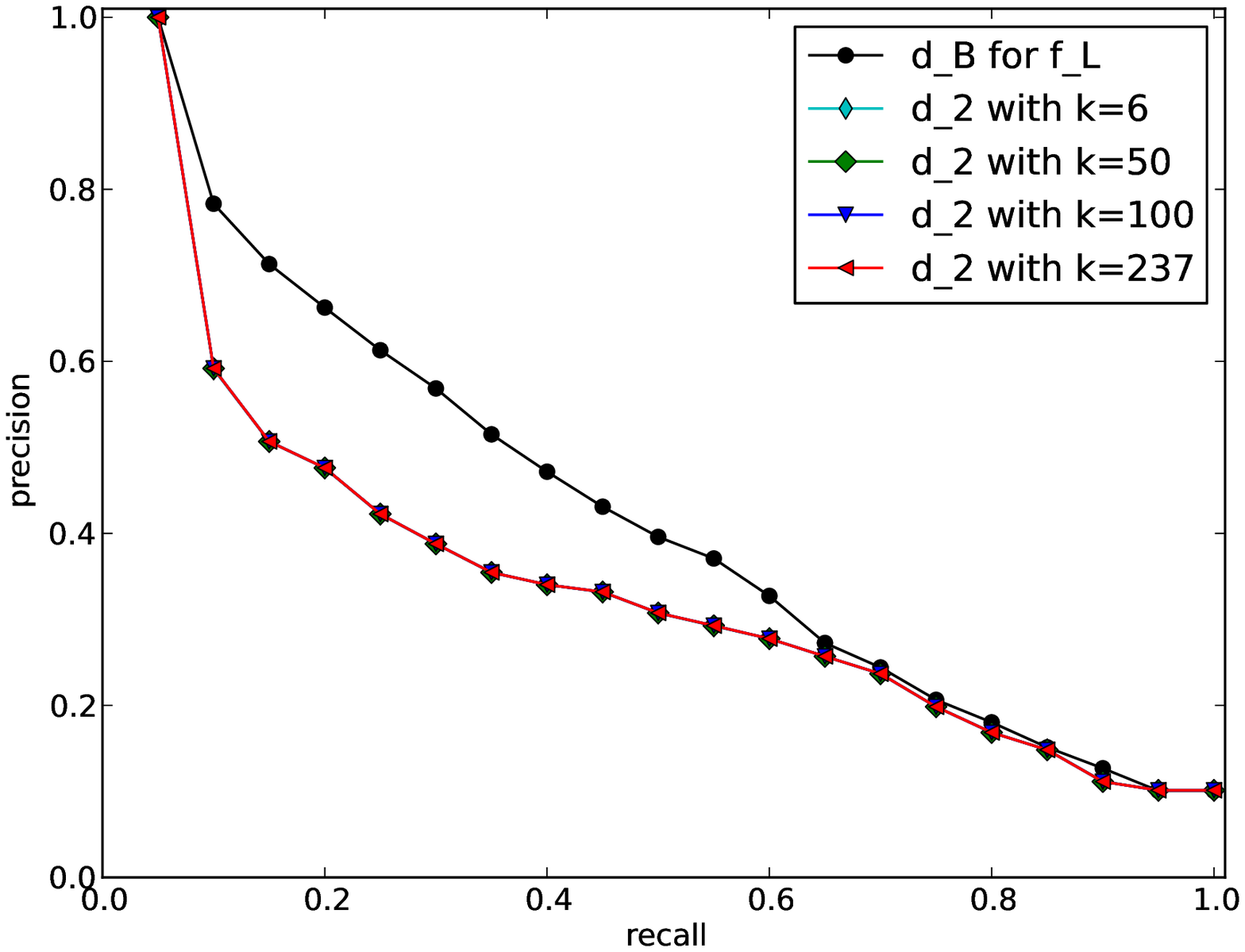}\!\!\!\!&\!\!\!
\includegraphics[height=3.2cm]{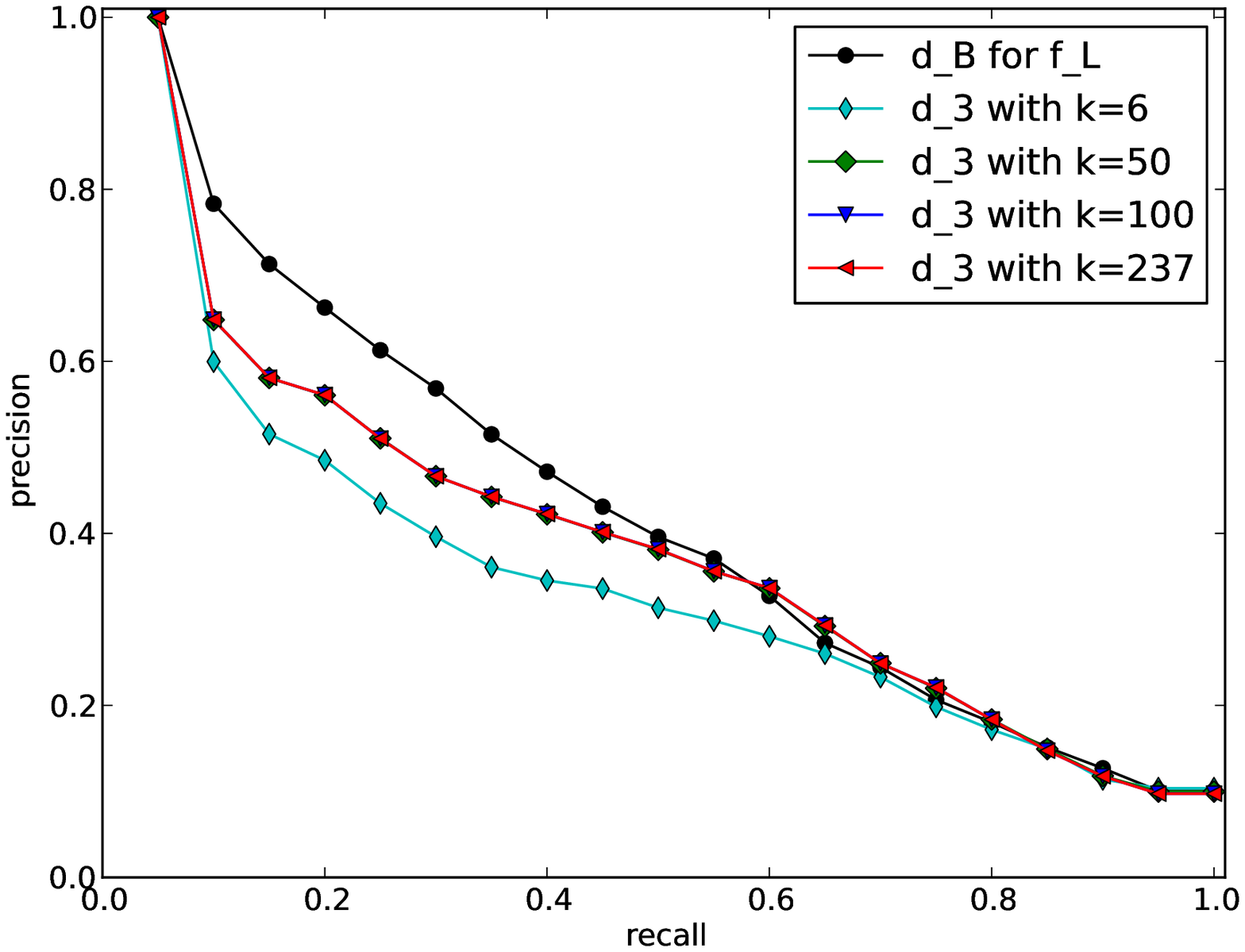}\\
\includegraphics[height=3.2cm]{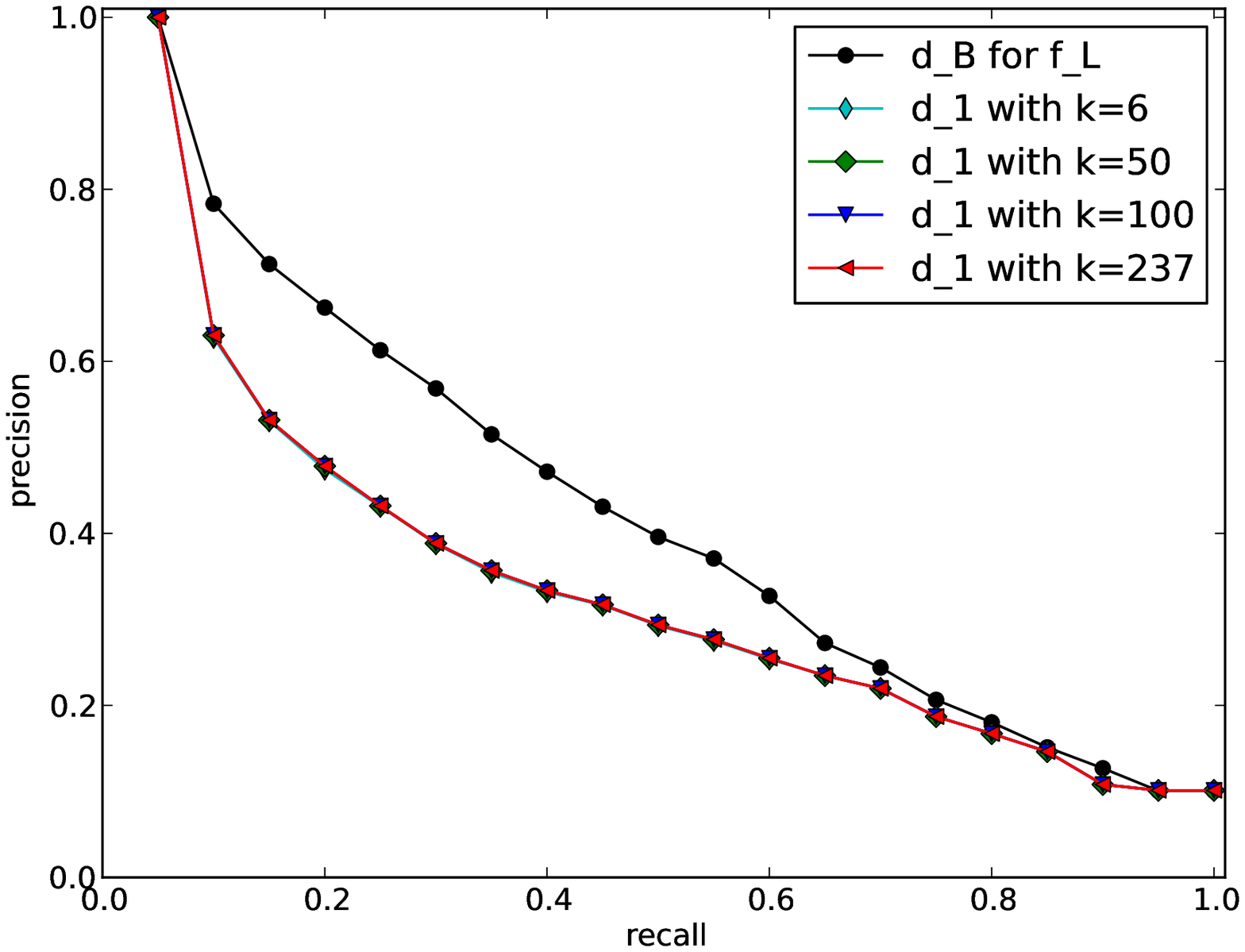}\!\!\!\!&\!\!\!
\includegraphics[height=3.2cm]{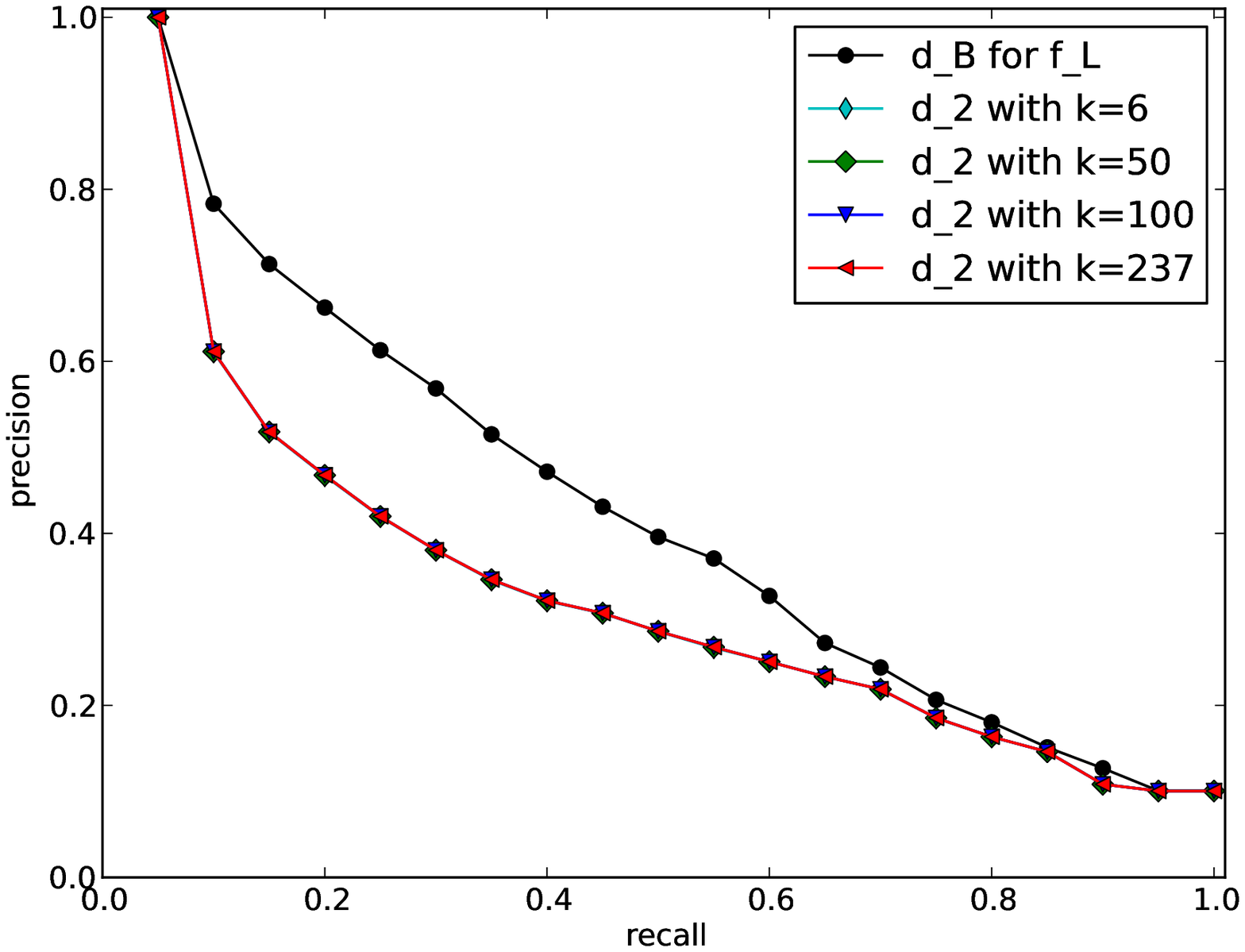}\!\!\!\!&\!\!\!
\includegraphics[height=3.2cm]{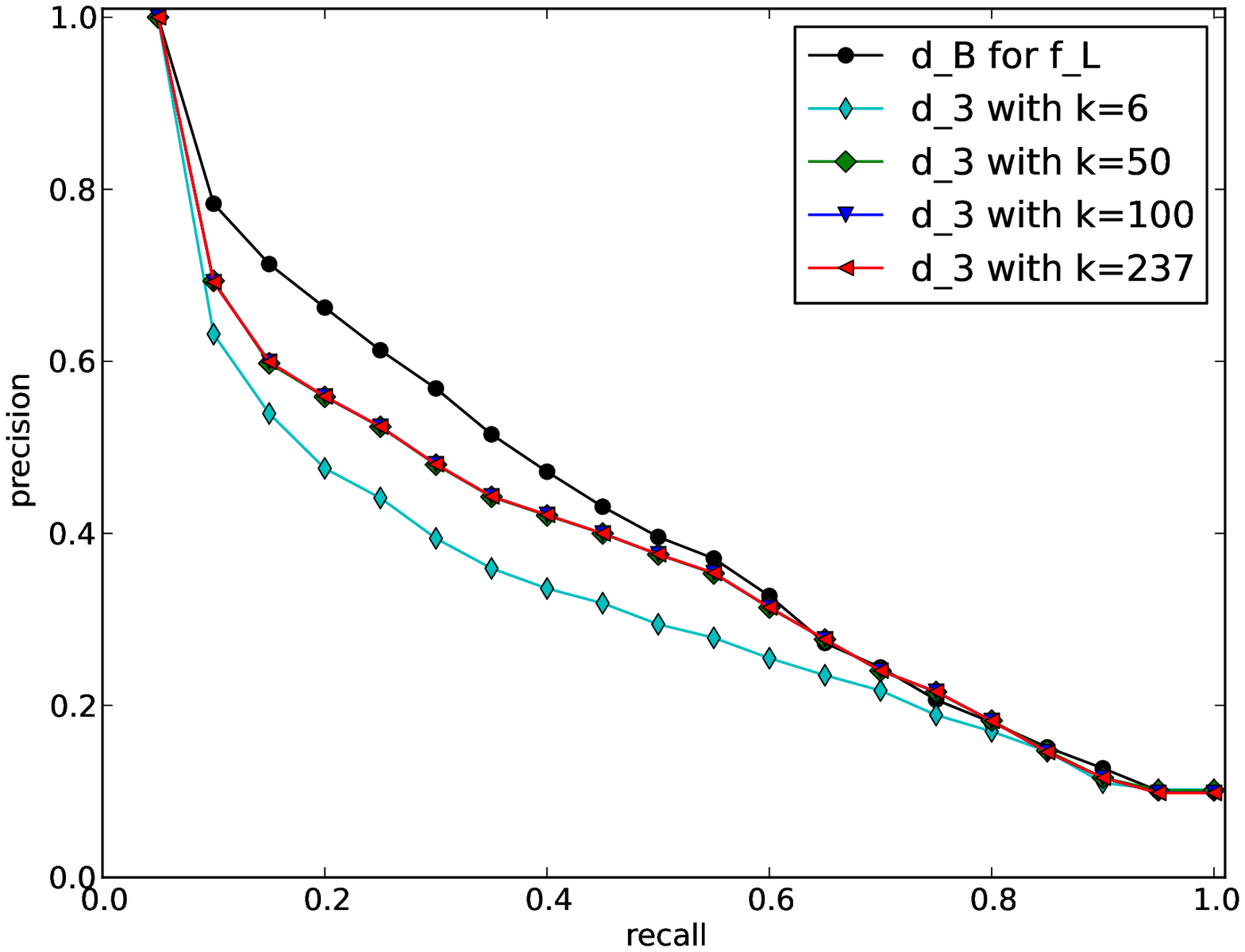}\\
\end{tabular}
\end{center}
\end{table}

\begin{table}[htbp]
\begin{center}
\caption{\footnotesize{PR graphs related to the filtering function $f_P$, when $0th$-persistence diagrams are compared directly through the bottleneck distance and in terms of the first $k$ components of the complex vectors obtained from the transformations $R$ (first row), $S$ (second row) and $T$ (third row) through the distances $d_1$ (first column), $d_2$ (second column) and $d_3$ (third column).}}\label{graphP}
\begin{tabular}{ccc}
\includegraphics[height=3.2cm]{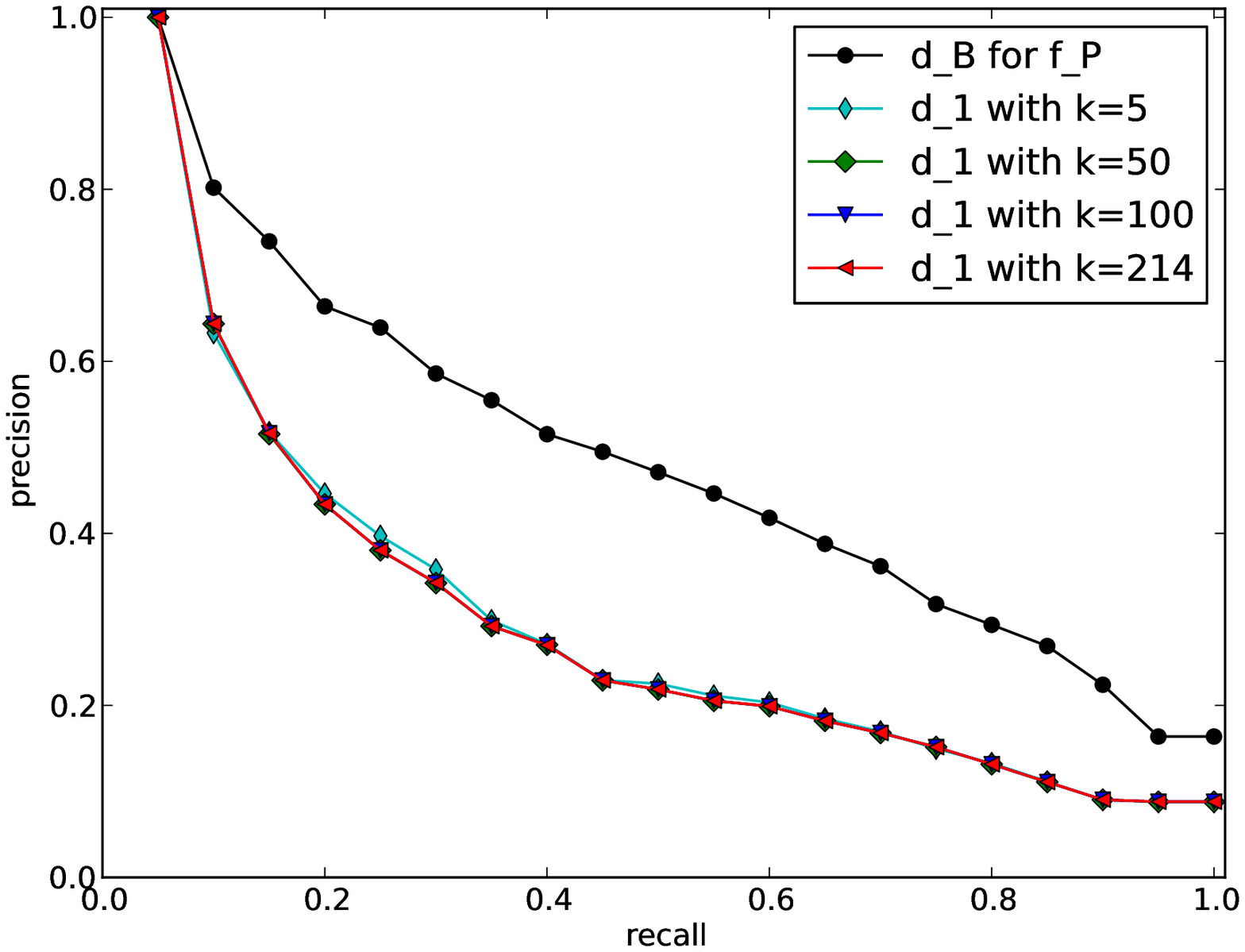}\!\!\!\!&\!\!\!
\includegraphics[height=3.2cm]{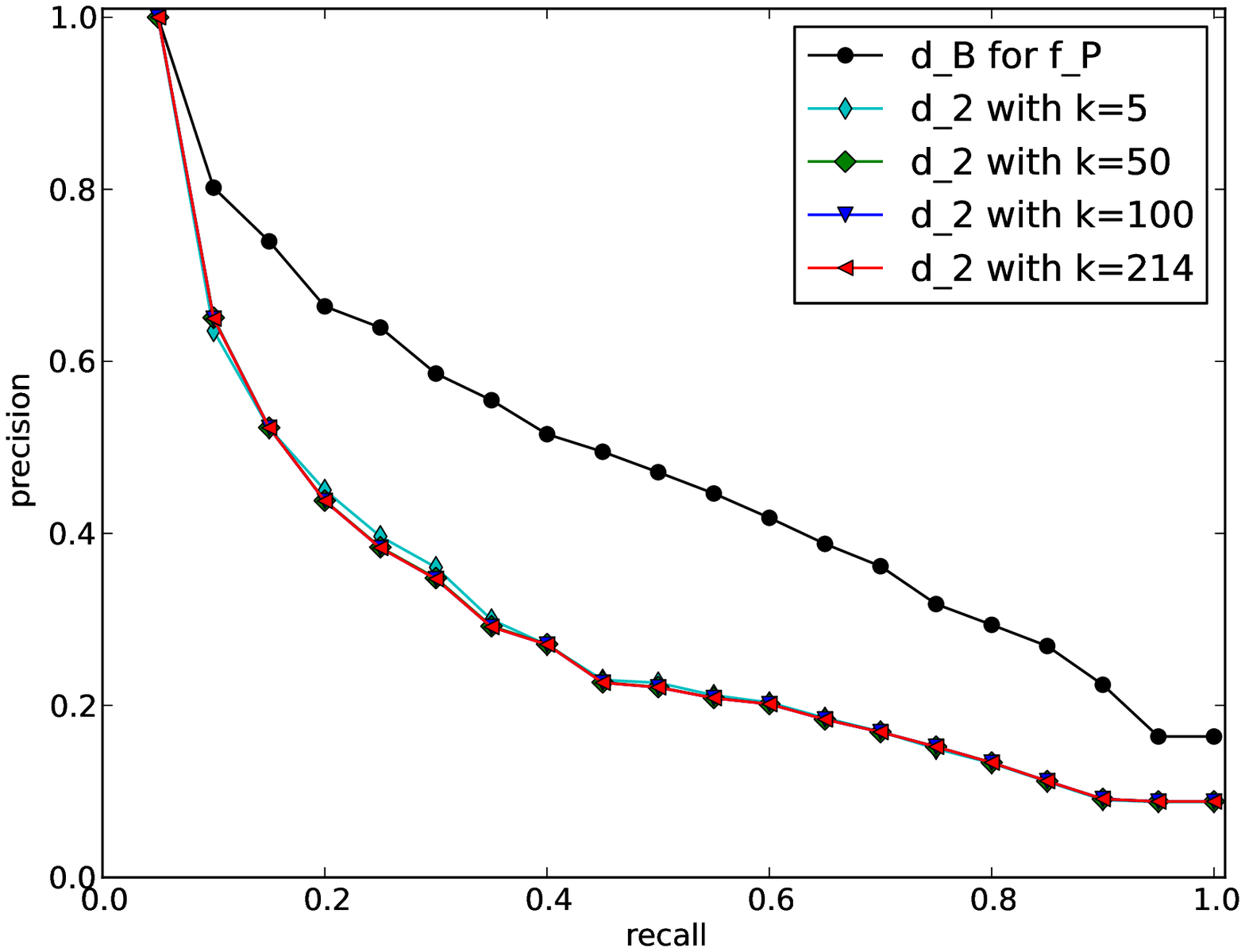}\!\!\!\!&\!\!\!
\includegraphics[height=3.2cm]{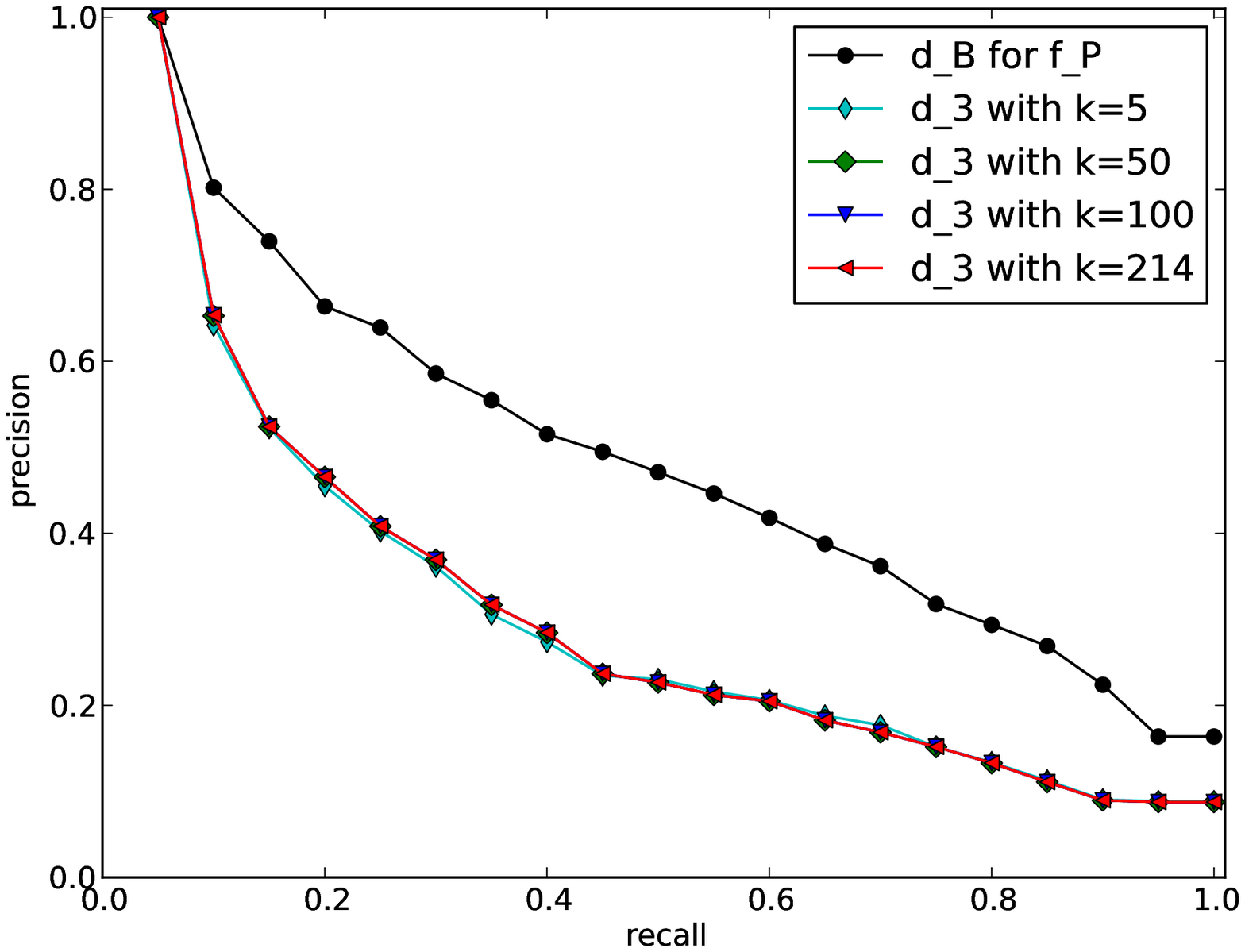}\\
\includegraphics[height=3.2cm]{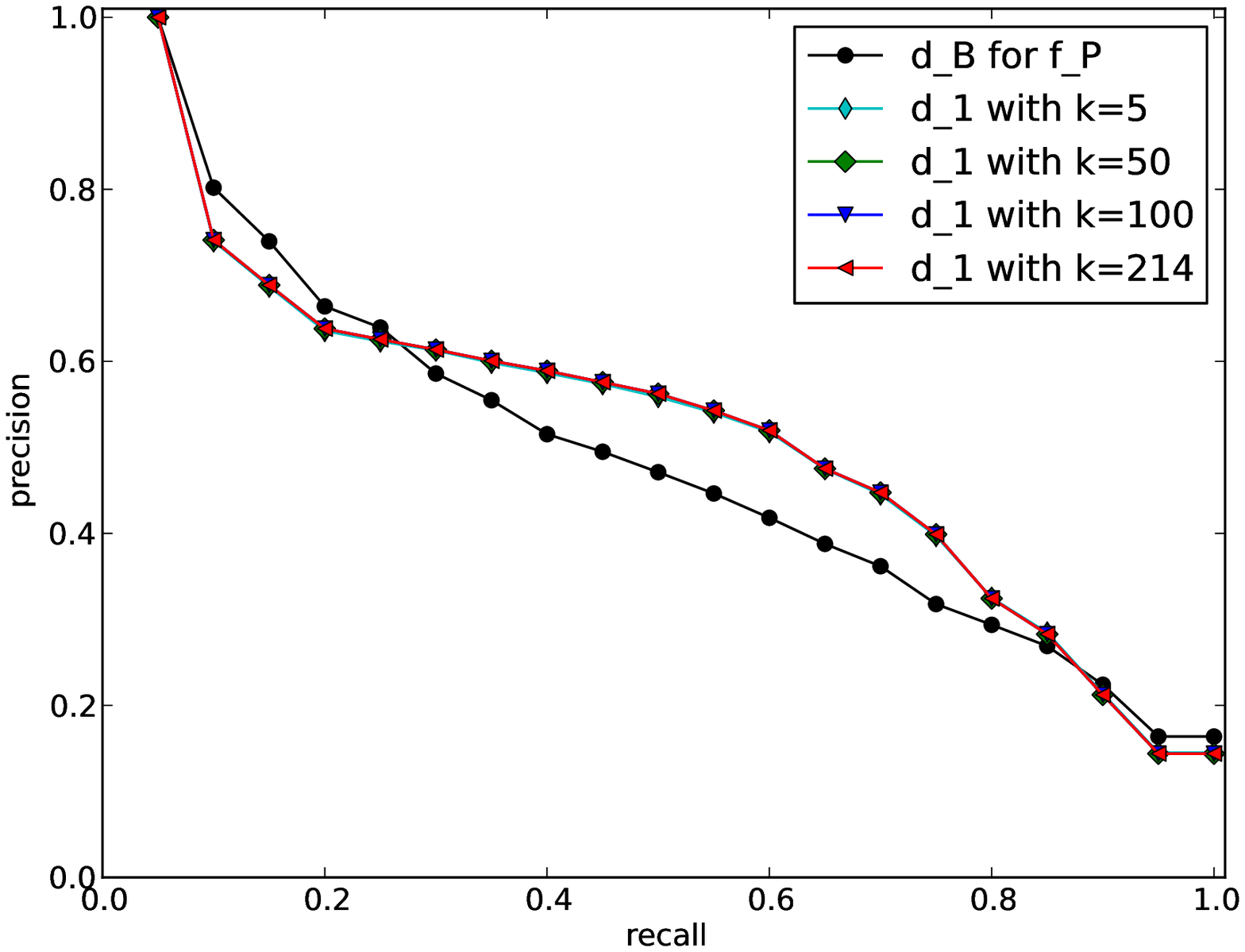}\!\!\!\!&\!\!\!
\includegraphics[height=3.2cm]{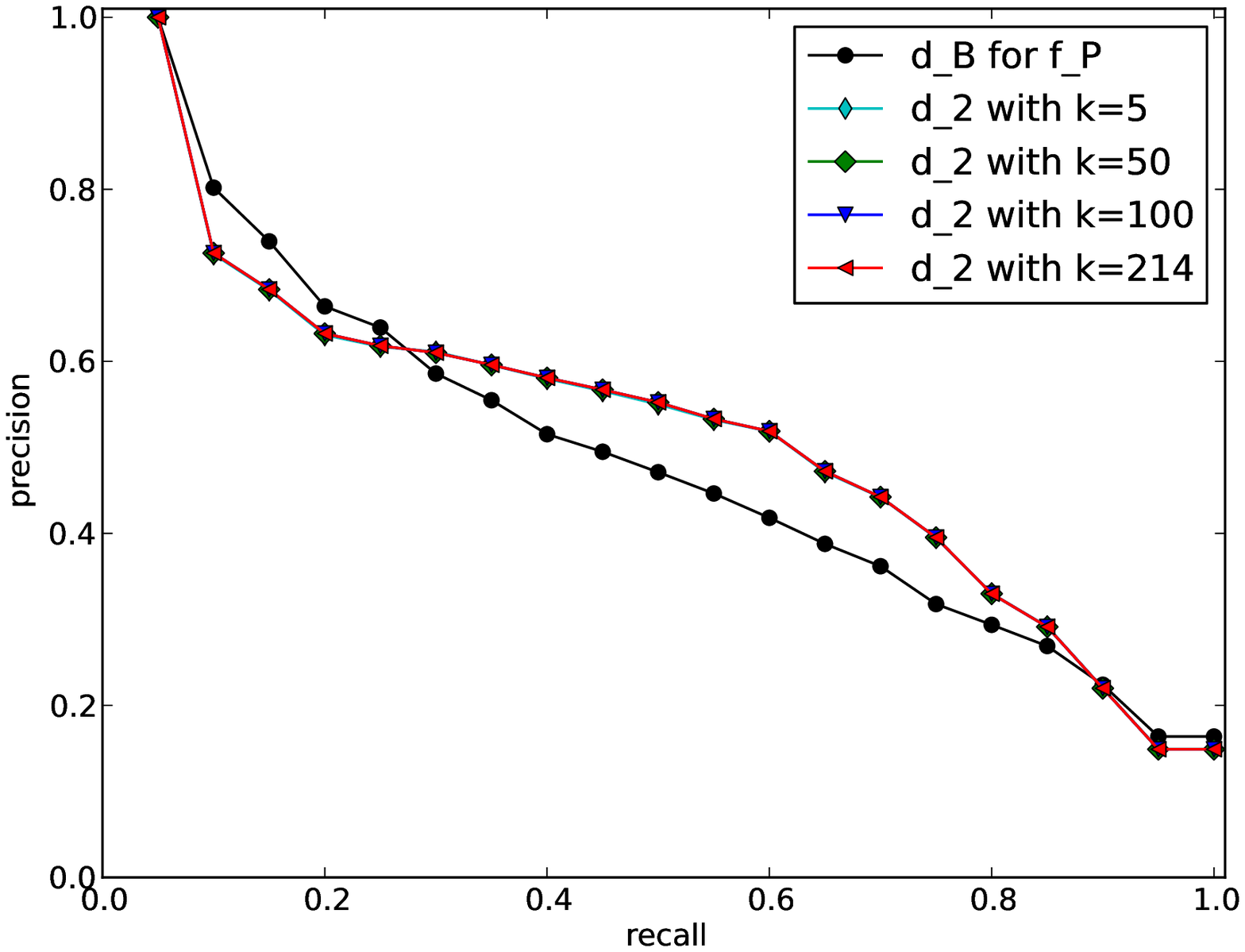}\!\!\!\!&\!\!\!
\includegraphics[height=3.2cm]{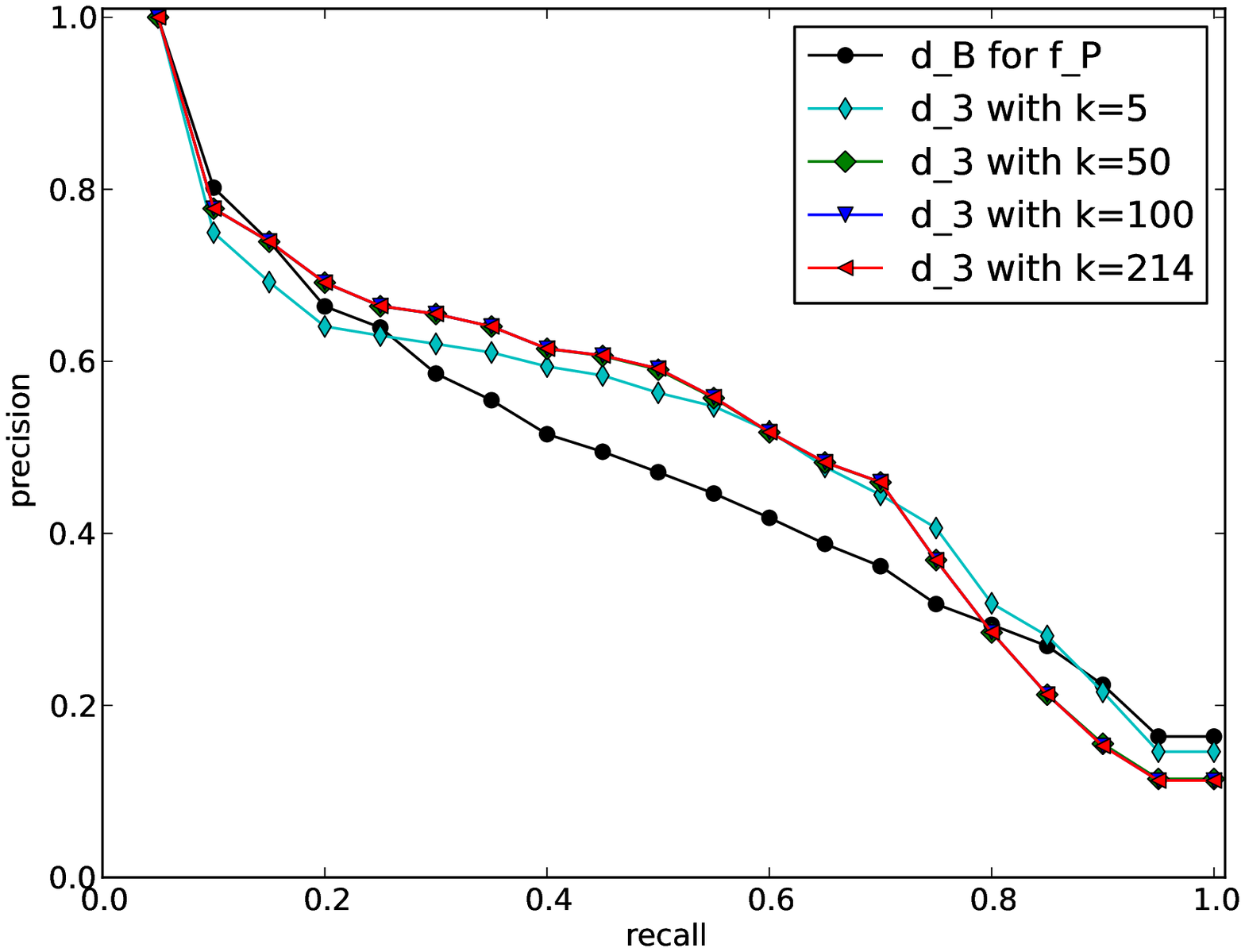}\\
\includegraphics[height=3.2cm]{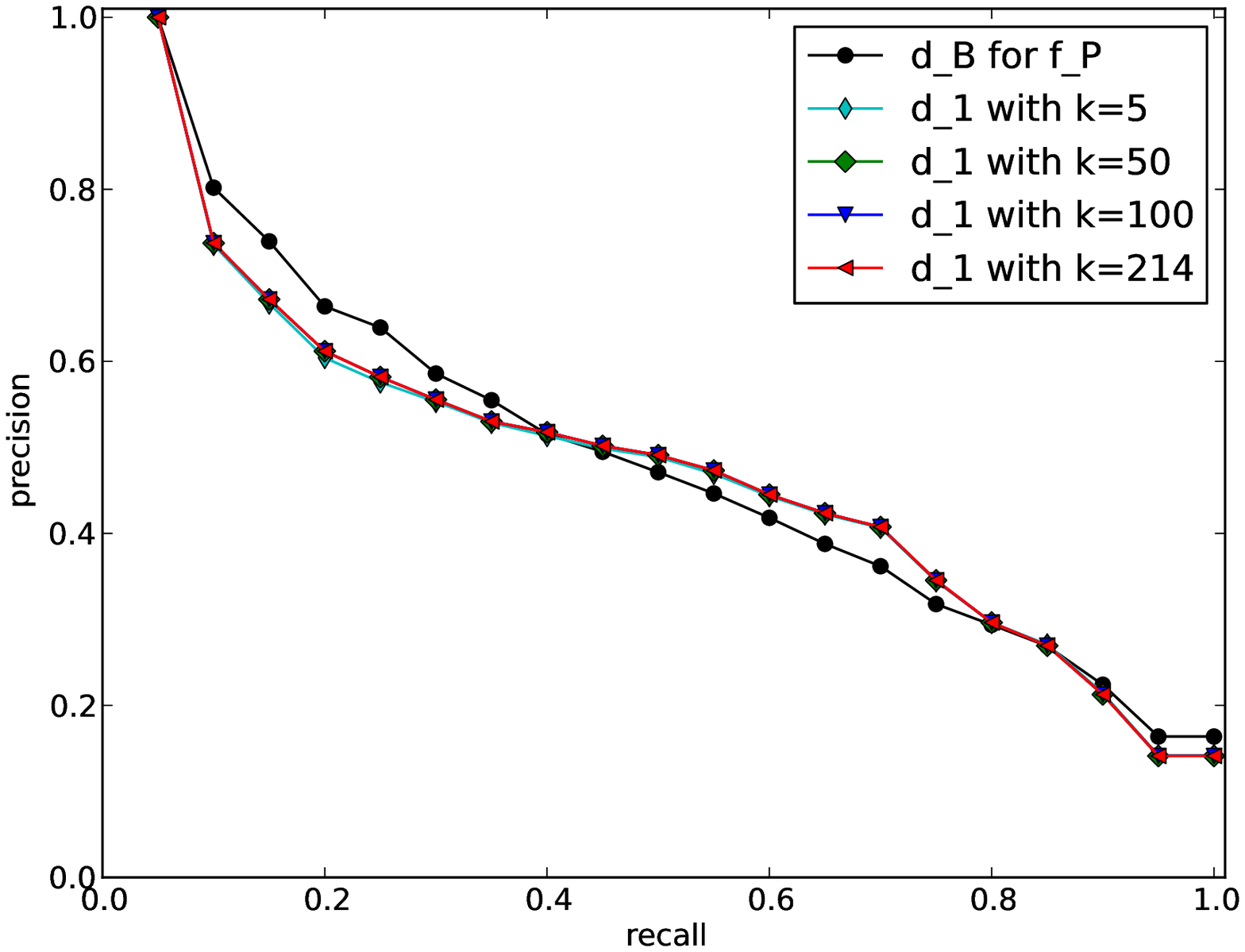}\!\!\!\!&\!\!\!
\includegraphics[height=3.2cm]{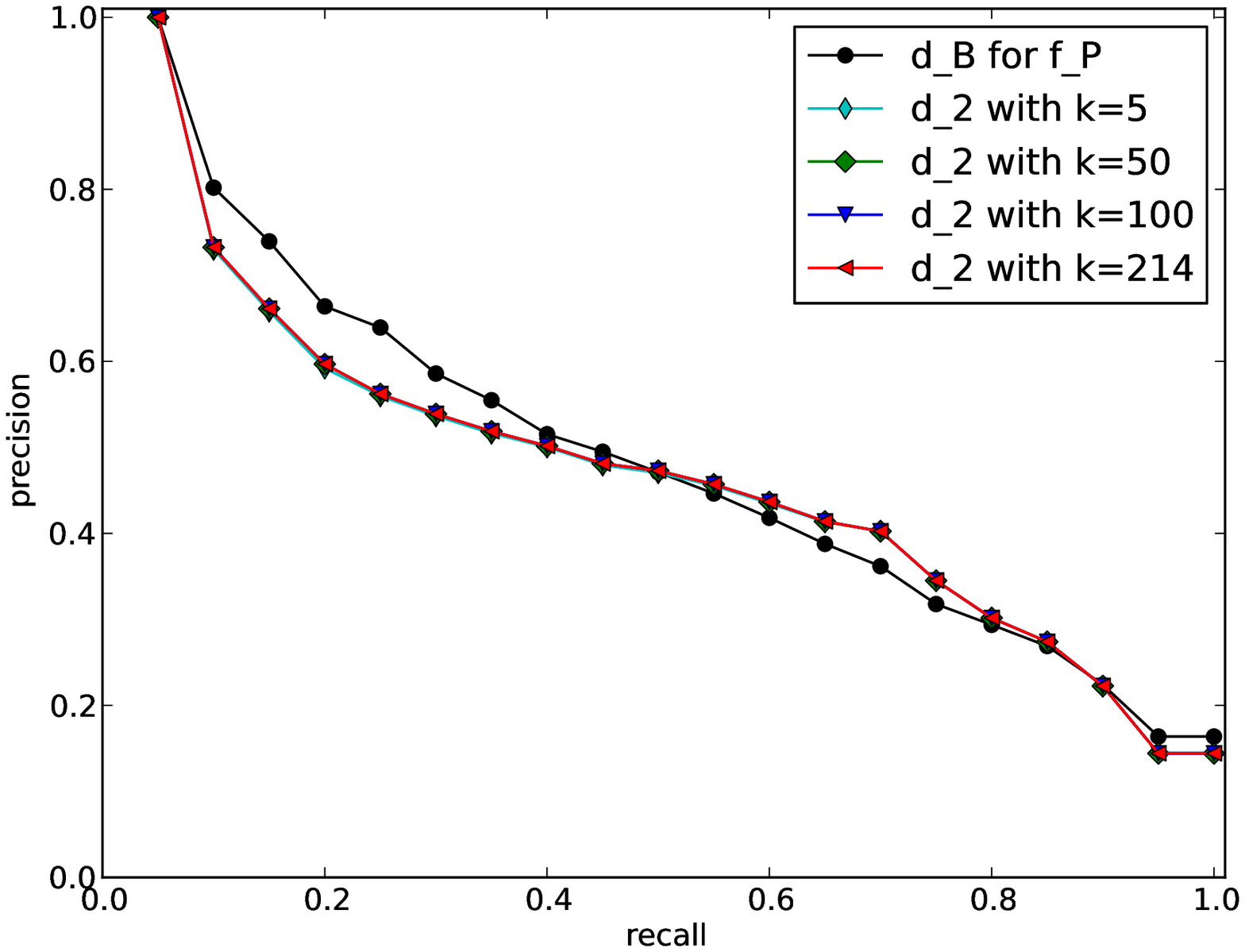}\!\!\!\!&\!\!\!
\includegraphics[height=3.2cm]{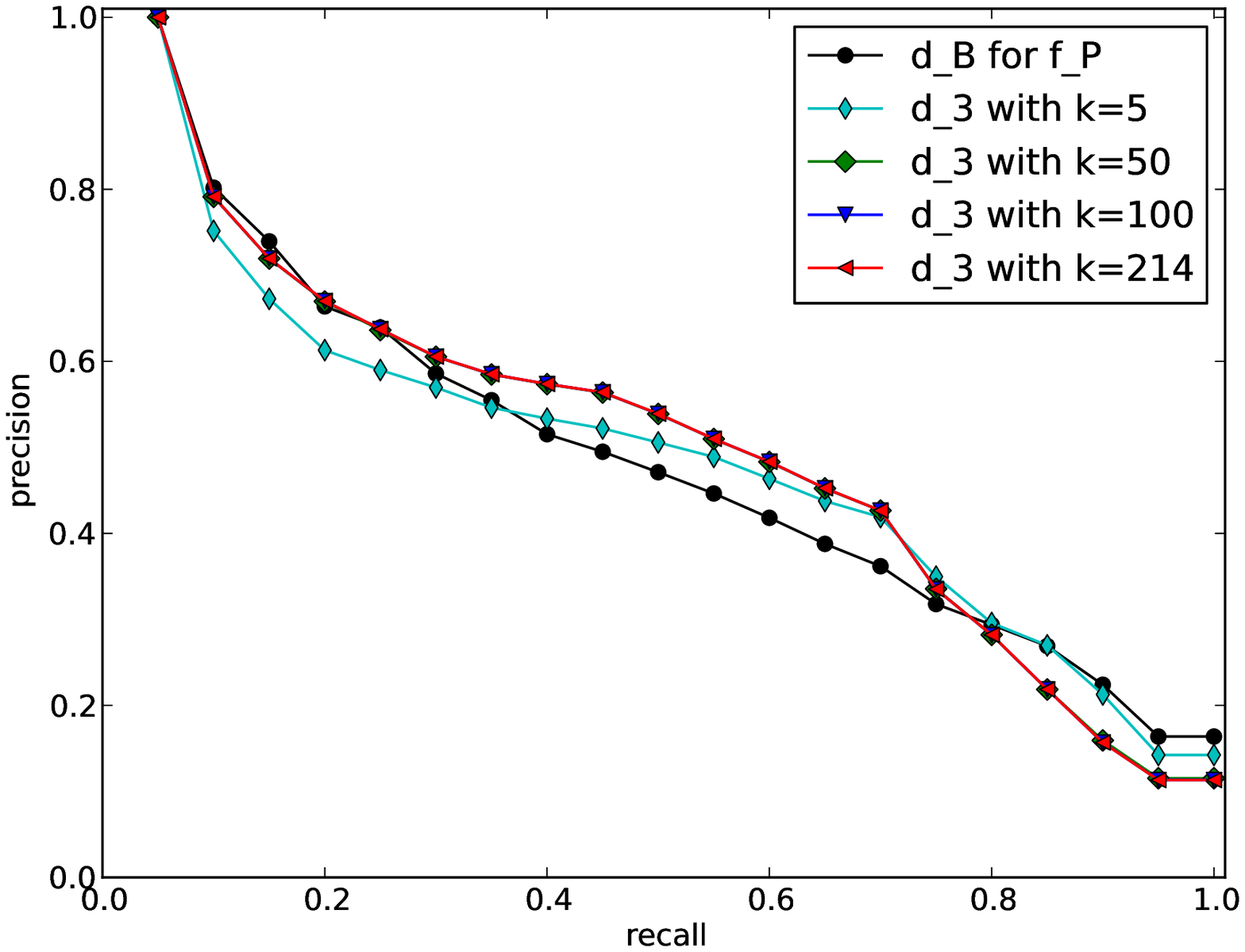}\\
\end{tabular}
\end{center}
\end{table}

What these plots aim to show is the  comparison of the
performances when the database classification is carried out
through the computation of the bottleneck distance $d_B$ between
persistence diagrams or the computation of the distances $d_1$,
$d_2$ and $d_3$ between the first $k$ components of the complex
vectors obtained through the transformations $R$, $S$ and $T$, for
different values of $k$  (see Section \ref{vectors} for the
definitions of $R$, $S$, $T$, $d_1$, $d_2$ and $d_3$). As it can
be easily observed, increasing the value of $k$ from the smallest
to the biggest number of proper points in the persistence diagrams
of our database, the PR graphs do not change so sensibly. This
means that the most important information of the persistence
diagram is contained in the first few vector components, the ones
corresponding to the coefficients of monomials with highest
degree. Moreover, we point out also that PR graphs related to
vectors which are induced by transformations warping the diagonal
$\Delta$ to a point (second and third rows in Tables \ref{graphL}
and \ref{graphP}) provide better results than by acting as the
identity (first row). This fact depends on the properties of
polynomial coefficients: Indeed roots corresponding to points of
persistence diagrams farther from the diagonal weigh more than
those closer to it. Hence, applying transformations $S$ and $T$
corresponds, in some sense, to providing points of a persistence
diagram with a weight that follows the paradigm of persistence:
The longer the lifespan of a homological class, the higher the
weight associated with the point having as coordinates the birth
and death dates of this class. This outcome is moreover
strengthened by the usage of the distance $d_3$ (third column in
Tables \ref{graphL} and \ref{graphP}) since it greatly enhances
the contribution of the first vector components to their
dissimilarity measure at the expense of the last.

Finally, note that the precision values at high recall --- i.e. by
retrieving a large number of objects --- are always fairly
comparable with the values relative to the bottleneck distance.
This assures us that complex vector comparison can act as a fast
and reliable preprocessing scheme for reducing the set of objects
to be fed to the generally more precise bottleneck distance.


{\small
\bibliographystyle{splncs03}
\bibliography{RefsICIAP}
}
\end{document}